\documentclass{ifacconf}

\usepackage{graphicx}
\usepackage{amsmath,amssymb}
\usepackage{natbib}  
\usepackage{xcolor}

\newtheorem{expl}{Example}
\newtheorem{ass}{Assumption}
\newcommand{\R}{{\mathbb R}}

\newcommand{\N}{{\mathbf N}}

\newcommand{\cL}{{\mathcal L}}

\newcommand{\cO}{{\mathcal O}}

\newcommand{\bE}{{\mathbf E}}

\newcommand{\one}{{\mathbf 1}}

\newcommand{\eps}{\epsilon}
\newcommand{\KL}{\it KL}

\newcommand{\wrt}{with respect to }

\begin{document}
\begin{frontmatter}

\title{Risk-neutral limit of adaptive importance sampling of random stopping times}


\author{Carsten Hartmann and Annika Jöster} 

\address{Brandenburgische Technische Universität Cottbus-Senftenberg, Institute of Mathematics, Konrad-Wachsmann-Allee 1, 03046 Cottbus, Germany (e-mail: carsten.hartmann@b-tu.de, annika.joester@b-tu.de)}

\begin{abstract}                
We discuss importance sampling of exit problems that involve unbounded stopping times; examples are mean first passage times, transition rates or committor probabilities in molecular dynamics. The naive application of variance minimization techniques can lead to pathologies here, including proposal measures that are not absolutely continuous to the reference measure or importance sampling estimators that formally have zero variance, but that produce infinitely long trajectories.  We illustrate these issues with simple examples and discuss a possible solution that is based on a risk-sensitive optimal control framework of importance sampling.
\end{abstract}

\end{frontmatter}

\section{Introduction}

The simulation of rare events is a challenge for Monte-Carlo simulations; see \cite{Asmussen2013}. To illustrate the typical difficulties, we consider $X_t\in\R^n$ solving  the stochastic differential equation (SDE)
\begin{equation}\label{ou}
    dX_t = -\theta X_tdt + \sqrt{2\eps}dW_t\,,\quad X_0=x\,,
\end{equation}
for some $\theta>0$ and small $\eps>0$, and suppose we are interested in the statistics of the exit of the process from the set  $B :=\{x\in\R^n\colon |x|<R\}$ for some $R>0$. We call  
\begin{equation}
    \tau_B := \inf\{t>0\colon X_t\notin B\}\,
\end{equation}
the first exit time of $B$. Quantities of interest to characterise the rare exit statistics are, e.g., 
the exit probability $P(\tau_B\le T)$ or the mean first exit time $\bE[\tau_B]$. For small $\eps$, the latter is exponentially large in $\eps$, i.e. 
\[
\lim_{\eps\to 0} \eps \log \bE[\tau_B] = \theta RL^2\,,
\]
independent of the initial condition $X_0=x$. Conversely, the exit probability is  exponentially small, which is consistent with the observation that the exit time asymptotically follows an exponential distribution; more precisely (see \cite{Zabczyk1985,IHP2019}):     
\[
\lim_{T\to\infty} \frac{1}{T}\log P(\tau_B>T) = -e^{-\frac{\theta R^2}{\eps}}\,. 
\]
As a consequence, the plain vanilla Monte Carlo estimator of the exit probability inevitably suffers from a relative error that is unbounded in $\eps$ (and diverges exponentially as $\eps\to 0$). Indeed, letting \[\tau_B^{(1)},\ldots,\tau_B^{(N)}\sim\tau_B\] denote $N$ independent copies of $\tau_B$, then
\begin{equation}
    \hat{p}_{N}(T) = \frac{1}{N}\sum_{i=1}^N \one_{\{\tau_B^{(i)}\le T\}}
\end{equation}
is an unbiased estimator of $P(\tau_B\le T)$, with relative error
\begin{equation}
    \delta(\eps):=\frac{\sqrt{{\rm Var}(\hat{p}_{N}(T))}}{\bE[\hat{p}_{N}(T)]}\asymp \frac{1}{\sqrt{NT}} e^{\frac{\theta L^2}{2\eps}}\,
\end{equation}
as $\eps\to 0$.

\section{Adaptive importance sampling}\label{sec:ais}

Our little introductory example explains why variance reduction methods for rare events are essential because the relative error can only be controlled through the variance when the estimator is unbiased. 
One such variance reduction method is importance sampling that changes the probability distribution of the  underlying  process by adding an extra drift to (\ref{ou}) that drives the dynamics towards the boundary of the set $B$, so that the exit from the set is no longer rare. By carefully choosing the extra drift, it is possible to decrease the variance of the estimator. Yet the approach may produce poor results or become inefficient for bad choices of the extra drift, since the likelihood ratio that appears as a multiplicative correction term in the corresponding importance sampling estimator is not bounded a priori.

Methods to adaptively choose the drift so as to minimise the estimator variance are known by the name of \emph{adaptive importance sampling (AIS)}. AIS methods that go back to the seminal works of \cite{DupuisWang2004,DupuisWang2007} are based on a (deterministic) control interpretation of the large deviations rate function associated with a rare event in the limit $\eps\to 0$; see \cite{Weare2012} for related work. The nonasymptotic variant  for $\eps>0$ that is based on a stochastic control framework using the Gibbs principle for the energy associated with the rare event has been suggested in \cite{Hartmann2012,zhang2014applications,HartmannEnt2017}. 

\subsection{Zero-variance AIS and the Feynman-Kac formula}

The holy grail of any AIS scheme are zero-variance estimators that can be constructed in certain cases. Let us briefly explain the idea for the case of the exit probability; the approach is based on an application of the celebrated Feynmac-Kac formula, following ideas described in \cite{awad2013,talayBook}. 

We consider a slight generalisation of the introductory example and assume that $X$ is an $n$-dimensional diffusion on $[0,\infty)$ governed by the SDE 
\begin{equation}\label{sde}
     dX_{t} = b(X_{t})dt + \sigma dW_{t}\,,\quad X_0=x\,,   
\end{equation}
with infinitesimal generator
\[
\cL = \frac{1}{2}\sigma\sigma^T\colon\nabla^2 + b\cdot\nabla \,.
\]
Here and in what follows we assume that $\sigma\in\R^{n\times n}$ is invertible and $b\colon\R^n\to\R^n$ is a smooth vector field that satisfies the usual (global) Lipschitz and growth conditions that guarantee that (\ref{sde}) has a unique strong solution. 

As before, we denote by $B\subset\R^n$ an open and bounded set, having smooth boundary $\partial B$ and call $\tau_B$ the corresponding first exit time. We moreover define the function $g(x)=\one_{\partial B}(x)$ and the stopping time $\tau:=\min\{\tau_B,T\}$ for some given $T\in(0,\infty)$, and we define the function 
\[
\psi(t,x)=\bE[g(X_{\tau})|X_t=x]\,.
\]
It can be readily seen that $\psi(t,x)$ equals the exit probability $P(\tau_B\le T|X_t=x)$. By the Feynman-Kac formula, it satisfies the linear parabolic PDE
\begin{equation}\label{FK}
\begin{aligned}
    	\left(\frac{\partial}{\partial t} + \cL\right) \psi(t,x) & = 0 \,, \quad (t,x)\in [0,T)\times B\\
     \psi(t,x) & = 1\,,\quad (t,x)\in [0,T)\times \partial B\\
        \psi(t,x) & = 0\,,\quad (t,x)\in \{T\}\times B\,.
\end{aligned}
\end{equation} 
Now suppose that we sample from another process 
\begin{equation}\label{controlledSDE}
    dX^u_t = \left(b(X^u_t)+\sigma u_t\right)dt + \sigma dW_t\,,\quad X^u_0=x\,,
\end{equation}
that depends on a control $u$ that satisfies suitable uniform integrability conditions that will be discussed a bit further below. Denoting by $\bE^Q[\cdot]$ the expectation \wrt $X^u$, we can trivially recast the exit probability as   
\begin{equation}\label{ISestim0}
    \bE[g(X_{\tau})] = \bE^Q\!\left[g(X_\tau)\frac{ \bE[g(X_{\tau})]}{g(X_{\tau})}\right]\,.
\end{equation}
By construction, the random variable 
\[
Z = g(X_\tau)\frac{ \bE[g(X_{\tau})]}{g(X_{\tau})}
\]
is constant, therefore its variance under the probability measure $Q$ is zero. (This property holds for all reasonable choices of $u$ as long as $Q=Q(u)$ is well-defined.) Clearly, equation (\ref{ISestim0}) is a tautology of very little practical use, but it suggests to find a representation of the control $u$, such that $Q$ and $P$ are mutually absolutely continuous, with likelihood ratio $L=dQ/dP$ given by 
\[
L= \frac{g(X_\tau)}{\bE[g(X_{\tau})]} = \exp(\log\psi(\tau,X_\tau) - \log\psi(0,X_0))\,,
\]
where the second equality follows from $\psi(\tau,X_\tau)=g(X_\tau)$ and $\psi(0,X_0)=\bE[g(X_{\tau})]$. A more convenient representation of $L$ in terms of the controlled process (\ref{controlledSDE}) is obtained by applying It\^o's formula to $Y_t=\log\psi(t,X_t)$ together with the Feynman-Kac formula (\ref{FK}), which yields  
\begin{equation}\label{LR}
L = \exp\!\left(\int_0^\tau v^*_t \cdot dW_t + \frac{1}{2}\int_0^\tau \left|v^*_t\right|^2 dt\right)
\end{equation}
where we used the shorthand 
\begin{equation}\label{controlFK}
    v^*_t = \sigma^T\nabla \log\psi(t,X_t)\,.
\end{equation}
By Girsanov's Theorem \cite[Thm. 8.6.8]{oeksendal2003} the measure $Q^{*}=Q(v^*)$ is generated by the feedback controlled SDE (\ref{controlledSDE}) with control $u=v^*$.  
By construction, the corresponding AIS estimator has zero variance: 
\begin{equation}\label{ISestim1}
    \bE[g(X_{\tau})] = g(X^*_\tau) \exp\!\left(-\frac{1}{2}\int_0^\tau \left|v^*_t\right|^2 dt-\int_0^\tau v^*_t \cdot dW_t \right)\,.
\end{equation}

\begin{rem}
    As the optimal control $v^*$ in (\ref{ISestim1}) depends on the quantity of interest, $\psi$, the zero-variance AIS scheme is not feasible. Nevertheless it is possible to reduce the variance by replacing $\psi$ by a suitable approximation. We will come back to this point in Section \ref{sec:num}.
\end{rem}

\subsection{Issues for unbounded stopping times}

The previous considerations readily extend to more general functionals of the process $X$. Specificlly, we consider 
	\[ 
	S(X) = \int_{0}^{\tau} f(X_{s})\,ds + g(X_{\tau})\,,
	\]
	for suitable functions $f,g\ge 0$ and an almost surely finite stopping time $\tau$. Its mean then  involves, e.g. 
	\begin{itemize}
		\item committor probabilities $P(\tau_B<\tau_A)=\bE[S(X)]$ for $f=0$, $g=\one_B$, $\tau=\min\{\tau_A,\tau_B\}$ where $A\subset\R^n$ is another subset that is disjoint from $B$
		\item mean first passage times $\bE[\tau_C]=\bE[S(X)]$ of a set $C\subset\R^n$ for $f=1$, $g=0$, and $\tau=\tau_C$ 
  
		\item hitting probabilities $P(\tau_B\le T)=\bE[S(X)]$ for $f=0$, $g(x)=\one_B(x)$ and $\tau=\min\{\tau_B,T\}$. 
	\end{itemize}

What distinguishes the first two examples from the third one is that in the latter case the length of the trajectories is at most $T$, whereas the trajectories in the first two cases can be arbitrarily long. As has been discussed in \cite{awad2013} and \cite{anum}, variance reduction (i.e.~reduction of the second moment) does not automatically lead to a reduction of the average length of the trajectories. This means that, while the AIS estimator under the  controlled dynamics (\ref{controlledSDE}) with control
\[
v^*_t = \sigma^T\nabla \log \bE[S(X)|X_t=x]
\]
enjoys a zero-variance property, this does not necessarily lead to an increase of the computational efficiency of AIS. We shall illustrate this aspect by revisiting the example from the introduction. 

\begin{expl}\label{ex:mfet}
    Consider (\ref{ou}) where for simplicity we set $\theta=0$. The mean first exit time (MFET) from the open ball $B$ of radius $R$ is the unique solution to the linear boundary value problem \cite[Thm.~9.1.1]{oeksendal2003}
    \begin{align*}
         -\cL\phi(x) & = 1\,,\quad x\in B \\
         \phi(x) & = 0\,,\quad x\in \partial B\,,
    \end{align*}
    with $\cL=\eps\Delta$. Its solution is easily found, exploiting the spherical symmetry of the problem: 
\[
\phi(x) := \bE[\tau_B|X_t=x] = \frac{R^2-|x|^2}{2n\eps}.
\]
\end{expl}
By the strong Markov property, the MFET is independent of the starting time $t$, and the AIS dynamics (\ref{controlledSDE}) with stationary control $v^*_t=\sqrt{2\eps}\nabla\log\phi(X_t)$ reads  
\begin{equation}\label{controlledBM}
	dX^*_{t} =  -\sqrt{2\eps}\frac{X^*_t}{R^2-|X^*_t|^2}dt + \sqrt{2\eps}dW_{t}\,,\quad X^*_{0}=x\,.
\end{equation}
Formally, the controlled dynamics realises a zero-variance AIS scheme. Yet, the dynamics is mean reverting towards the origin where the drift becomes singular at the boundary of $B$, preventing the trajectories from ever reaching the boundary $\partial B$.  As a consequence, the corresponding zero-variance probability measure $Q^*$ has the property 
\[
Q^*(\tau_B=\infty)=1\,.
\]
Not only does this imply that $P$ and $Q^*$ are mutually singular, since $P(\tau_B=\infty)=0$, but $Q^*(\tau_B=\infty)=1$, and $Q^*(\tau_B<\infty)=0$, but $P(\tau_B<\infty)=1$; but it also implies that AIS with the optimal control $u^*$ would require infinite computational time, which renders AIS useless from a computational perspective.

\section{Control of moments}

The previous example illustrates an essential weakness of importance sampling as a method to control the variance of path properties when the length of the paths is not bounded. The observation that the $Q$-variance is finite or even zero while the first moment diverges is not a contradiction, because AIS controls only control the second moment of $\tau L^{-1}$ under $Q$, which implies a bound on the first moment of $\tau L^{-1}$, but not on the first moment of $\tau$.  

The loss of control over the first moment is less severe when the stopping time $\tau$ is bounded, since in this case the maximum length of the sampled trajectories and hence the average simulation time are bounded.   
Therefore, we extend our AIS wish list for problems involving unbounded stopping times in that we aim at finding a change of measure from $P$ to $Q$ that \emph{both} reduces the estimator variance \emph{and} the average length of the controlled trajectories. 

\begin{rem}
    For the problem considered in Example \ref{ex:mfet}, a simple fix to avoid the singularity at the boundary of the target set is to replace $\bE[\tau_B]$ by $\bE[\tau_B+c]$ where $c>0$ is some regularisation parameter that bounds the quantity of interest away from zero and that can be subtracted from the AIS estimator to get unbiased minimum-variance estimate of the MFET. Since this form of regularisation does still \emph{increase} the average trajectory length, rather than \emph{decreasing} it, it can be considered only a partial solution to the problem. 
\end{rem}

\subsection{A certainty equivalence principle}

We can control several moments at once by introducing a risk-sensitive formulation of the problem \cite[Sec.~1]{whittle2002}; instead of $\bE[S(X)]$, we consider the certainty-equivalent expectation
\begin{equation}\label{cert}
    \gamma = \varphi^{-1}(\bE[\varphi(S(X))])
\end{equation}
where $\varphi$ is a strictly convex (strictly increasing or decreasing) function with inverse $\varphi^{-1}$. 

Two notable special cases of (\ref{cert}) are

	\begin{itemize}
		\item $\varphi(s)=|s|^p$ for $p > 1$, with the property
		\[
		\left(\bE\big[(S(X))^p\big]\right)^{1/p} \ge \bE\big[S(X)\big]
		\]
		
		\item $\varphi(s)=e^{-\alpha s}$ for $\alpha>0$, with the property
		\[
		-\alpha^{-1}\log\bE\big[e^{-\alpha S(X)}\big] \le  \bE\big[S(X)\big]\,.
		\]
	\end{itemize}
Since $\varphi$ is \emph{strictly} convex, equality holds iff $S$ is a.s.~constant, and we can use this fact as a characterisation of a change of measure that nullifies the variance. (A random variable is constant iff its variance is zero.) 

We focus on the second case that allows us to simultaneously control all moments of $S(X)$ at once, due to its connection with the log moment or cumulant generating function (CGF).  To this end, we remind the reader of the definition of the CGF of a random variable $Y$ where we consider a scaled version of the traditional CGF, also known as free energy (cf.~\cite[Ch.~4.1]{kendall}):
\begin{equation}\label{cgf}
        \gamma\colon\R\to (-\infty,\infty]\,,\; \alpha\mapsto -\frac{1}{\alpha}\log\bE\big[e^{-\alpha Y}\big]\,. 
\end{equation}
Note that the CGF encodes information about the mean and the variance of the random variable $Y$, assuming they exit, in that for sufficiently small $\alpha>0$: 
\[
\gamma(\alpha) = \bE[Y] - \frac{\alpha}{2}{\rm Var}(Y) + \cO(\alpha^2)\,.
\]
Moreover,  if $Y\ge 0$ then $\gamma(\alpha)$ is finite for any $\alpha>0$, and Jensen's inequality implies that  
\begin{equation}
    \gamma(\alpha) \le \bE^Q[Y] + \frac{1}{\alpha}\KL(Q,P)
\end{equation}
for all $Q$ absolutely continuous \wrt $P$ (symbolically: $Q\ll P$) where
\begin{equation}
    \KL(Q,P) = \begin{cases}
        \displaystyle \int\log\left(\frac{dQ}{dP}\right)dQ &\textrm{if }\; Q\ll P\\
        \infty & \textrm{otherwise}
    \end{cases}
\end{equation}
denotes the relative entropy (or: Kullback-Leibler divergence) between $Q$ and $P$. By the strict convexity of the exponential function $\varphi(s)=\exp(-\alpha s)$, equality is attained iff the random variable $Y+\alpha^{-1}\log\big(\frac{dQ}{dP}\big)$ is $Q$-a.s. constant. 

For finite-time SDE path functionals $Y=S(X)$, the CGF admits a variational characterisation that goes back to \cite{daipra1996,Boue1998}. The following result extends these results to AIS of path functionals that involve random stopping times:
\begin{thm}[\cite{HartmannEnt2017}]\label{thm:duality}
Assuming sufficient regularity of the coefficients  $f,g,b$ and an a.s. finite stopping time $\tau$, the CGF is the value function of the following optimal control problem: minimise
\begin{equation}\label{cost}
J(u) = \bE\!\left[ \int_{0}^{\tau}\left(f(X^u_{s}) + \frac{1}{2\alpha} |u_{s}|^{2}\right)ds + g(X^u_\tau)\right]
\end{equation}
subject to  $X^u$ being the solution of the controlled SDE (\ref{controlledSDE}) with initial conditions $X^{u}_{t}=x$. 
That is, $V(t,x)=\min_u J(u)$. The minimiser $u^*$ is unique and given by the feedback law 
\[u^*_t = -\alpha\sigma^T\nabla  V(t,X_t^*)\,.\] 
Moreover, with probability one, 
\begin{equation}\label{CGFais}
\gamma(\alpha) = S(X^*) + \frac{1}{\alpha}\int_0^\tau u^*_t \cdot dW_t + \frac{1}{2\alpha}\int_0^\tau \left|u^*_t\right|^2 dt\,.
\end{equation}
where $\gamma(\alpha):=\gamma(\alpha;t,x)$ denotes the CGF as a function of the initial data $(t,x)$ for any $\alpha>0$. The latter is equal to the value function, in other words: 
\[V(t,x)=\gamma(\alpha;t,x)\,.\]
\end{thm}

Theorem \ref{thm:duality} is a nonlinear generalisation of the AIS scheme from Section \ref{sec:ais}, in that the quantity of interest is now the ``nonlinear'' expectation $-\alpha^{-1}\log\bE[\exp(-\alpha S(X))]$ of $S(X)$ instead of $\bE[S(X)]$; the corresponding nonlinear Feynmac-Kac formula has the form of a Hamilton-Jacobi equation for the value function. 

A more common, but equivalent formulation of the identity (\ref{CGFais}) is obtained by switching from the CGF to the moment-generating function (MGF):
\begin{equation}\label{MGFais}
e^{-\alpha \gamma(\alpha)} = \bE^{Q^*}\!\left[e^{-\alpha S(X^*)}\frac{dP}{dQ^*}\right] = \e^{-\alpha S(X^*)}\frac{dP}{dQ^*}\,,
\end{equation}
with (inverse) likelihood ratio
\begin{equation}
    \frac{dP}{dQ^*} = \exp\left(-\frac{1}{\alpha}\int_0^\tau u^*_t \cdot dW_t - \frac{1}{2\alpha}\int_0^\tau \left|u^*_t\right|^2 dt\right)\,.
\end{equation}

\begin{expl}
With regard to our wish list, it is instructive to consider the special case $f=1$, $g=0$, and $\tau=\tau_B$ being the first exit time from $B$. In this case
\begin{equation}\label{dv2}
-\frac{1}{\alpha}\log\bE\!\left[e^{-\alpha\tau}\right] = \min_u\bE\!\left[ \tau^u + \frac{1}{2\alpha} |u_{s}|^{2}ds\right],
\end{equation}
where we denoted the stopping time on the right hand side by $\tau^u$ to emphasize its dependence on the control $u$. This shows that the control has two effects: it minimises the estimator variance, but it also reduces the controlled MFET (i.e.~the average simulation time per realisation). 
 
\end{expl}

Clearly, as before, a computationally feasible AIS scheme for the CGF or the MGF requires a suitable approximation of the value function (i.e.~the optimal control). See \cite{nuesken2021,suboptimalIS} regarding some algorithmic aspects.

\subsection{Extracting the mean}

Even though the CGF or the MGF encode the statistics of $S(X)$ in terms of all its cumulants or moments, we are normally interested in computing single moments, say, the mean $\bE[S(X)]$ rather than its CGF or MGF. Clearly, 
\[
- \lim_{\alpha\searrow 0}\frac{1}{\alpha}\log\bE\!\left[e^{-\alpha S(X)}\right] = \bE[S(X)]\,,
\]
but for small $\alpha$, the control in (\ref{cost}) becomes heavily penalised and vanishes as $\alpha\to 0$, which implies that AIS for small values of $\alpha$ can no longer lead to a reduction of the average simulation time. 

Intriguingly, even though the optimal control vanishes in the limit $\alpha\to 0$, the AIS estimator based on (\ref{CGFais}) has a nontrivial limit, since the zero-variance property in (\ref{CGFais}) holds for all $\alpha>0$ and the left hand side in (\ref{CGFais}) converges to $\bE[S(X)]$.  The next theorem formalises this observation: 

\begin{thm}\label{thm:is2cv}
    Let $L^*=\frac{dQ^*}{dP}$ be the likelihood ratio associated with the change of measure from the reference measure $P$ to the zero-variance measure $Q^*=Q(u^*)$: 
    \begin{equation}
    L^* = \exp\left(\frac{1}{\alpha}\int_0^\tau u^*_t \cdot dW_t + \frac{1}{2\alpha}\int_0^\tau \left|u^*_t\right|^2 dt\right)\,.
\end{equation}
Then, with probability one, 
\begin{equation}
    - \lim_{\alpha\searrow 0}\frac{1}{\alpha}\log\bE^{Q^*}\!\left[e^{-\alpha S(X)-\log L^*}\right] = S(X) - M_\tau(X),
\end{equation}
where 
\begin{equation}
    M_t(X) = \int_0^t \sigma^T\nabla\psi(s,X_s)\cdot dW_s
\end{equation}
is a martingale with the property 
\[
\bE[M_\tau(X)] = 0\,,\quad {\rm Var}(S(X) + M_\tau(X)) = 0\,.
\]
\end{thm}

The sketch of proof of the theorem is deferred to the Appendix; here we will briefly on the relevance of the theorem, which is a statement of the convergence of the CGF \emph{after} doing the change of measure from $P$ to $Q^*$ and which entails a convergence statement of the optimal control $u^*$ as $\alpha\to 0$, for rare event simulation: 

\begin{enumerate}
    \item Even though the CGF of $S(X)$ that agrees with the value function $V$ converges to the mean of $S(X)$, i.e.~the function $\psi$ associated with the linear Feynman-Kac formula (\ref{FK}), the scaled control 
    \[\frac{u^*}{\alpha}\propto\nabla V\] 
    does \emph{not} converge to the control 
    \[v^*\propto\nabla\log\psi\]
    that realises the zero variance change of measure for $\psi$, but to a control (variate) proportional to $\nabla\psi$.
    \item Since, without scaling, $u^*$ converges to 0 as $\alpha\to 0$, the optimally controlled dynamics $X^*$ for $\alpha>0$ becomes uncontrolled in the limit. As a consequence, the mean simulation time per realisation is no longer reduced. 
    \item Nevertheless, we obtain a proper zero-variance estimator for the mean that does not suffer from the issue of diverging first moment. Yet, there is no change of measure, because the dynamics is not altered. The variance reduction is due to the control variate term $M_\tau$ that is anticorrelated with $S(X)$ and therefore annihilates the variance of the estimator.      
\end{enumerate}

\section{Numerical illustration}\label{sec:num}

We briefly illustrate our findings with the standard toy example (\ref{ou}). To this end, we compute the MFET of a linear SDE of dimension $n=100$ from a ball of radius $R=10$, comparing standard Monte Carlo (MC), control variate (CoV), and a perturbed control variate (PCoV). 

\subsection{100-dimensional Brownian motion}

To begin with, we set $\theta=0$. The exact solution of the MFET in this case can be computed from solving the corresponding linear boundary value problem exploiting the spherical symmetry of the problem (cf.~Example \ref{ex:mfet}):   
\begin{equation}\label{mfet100}
\phi(x) = \bE\!\left[\tau|X_0=x\right] =\frac{100 - |x|^2}{200\eps}     
\end{equation}

\begin{figure}
    \centering
    \includegraphics[width=0.45\textwidth]{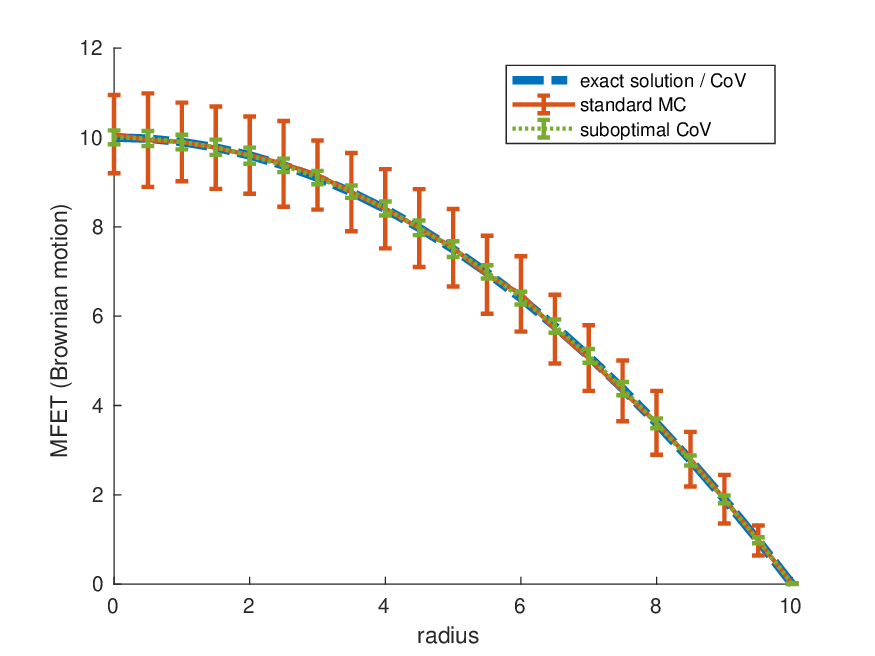}
    \caption{MFET estimates and their 95\% confidence intervals: exact solution that agrees with the zero-variance CoV estimate (dashed blue ), standard MC (solid red) and PCoV with perturbation $\delta\sin(x)$ for $\delta=0.25$ (dotted green). Simulations were done using the Euler-Maruyama scheme with steps size $\Delta t=10^{-4}$.}
    \label{fig:mfet}
\end{figure}

Here we have set $\eps=0.05$, which implies that the MFET ranges between 10 for a radius $|x|=0$ and 0 for $|x|=10$. Figure \ref{fig:mfet} shows the simulation results for $N=10$ independent realisations where the 95\% confidence intervals of the estimators have been computed from averaging over $M=100$ independent runs. As perturbation, we have added a component-wise $O(1)$ perturbation of the form $\delta\sin(x)$ with $\delta\in[0.25,1]$. Figure \ref{fig:mfet} only shows the PCoV result for $\delta=0.25$, where the relative error for larger values of $\delta$ remains within the range of PCoV and MC, but it never comes close to or even goes beyond the MC error.

\subsection{100-dimensional Ornstein-Uhlenbeck process}

\begin{figure}
    \centering
    \includegraphics[width=0.45\textwidth]{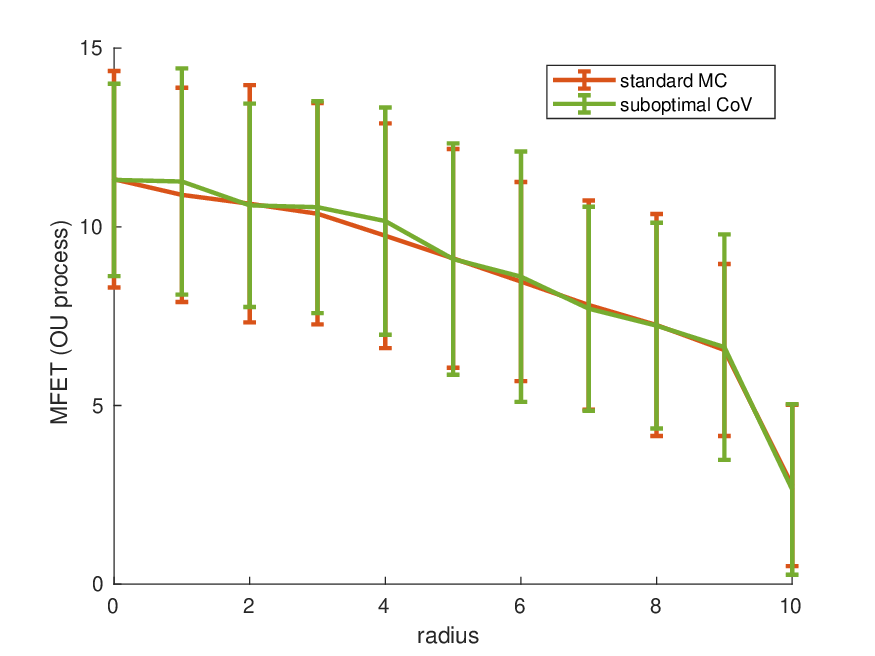}
    \caption{MFET estimates and their 95\% confidence intervals: standard MC (red) and CoV with a suboptimal CoV based on the driftless control variate based on (\ref{mfet100}). Simulations were done using the Euler-Maruyama scheme with steps size $\Delta t=10^{-3}$.}
    \label{fig:mfetOU}
\end{figure}

To further study the insensitivity of the CoV estimator to a bad approximation of the control variate we consider a 100-dimensional Ornstein-Uhlenbeck (OU) process of the form (\ref{ou}) where we replace the scalar $\theta>0$ by the symmetric positive definite drift matrix
\[
	\Theta = \begin{pmatrix}
		2 & -1 & 0 & \ldots & 0 \\ -1 &  2  & -1 & \ldots & 0\\ 0 & \ddots & \ddots & \ddots & 0\\ \vdots & \vdots & \ddots & \ddots  & -1 \\ 0 & \ldots & 0 & -1 & 2  
	\end{pmatrix}\in\R^{100\times 100}\,.
\]
Here we set again $\eps=0.05$. Explicit formulas for the MFET are only available in the spherically symmetric case (e.g.~\cite{sphericalOU}) or in the large deviations regime (e.g.~\cite{Zabczyk1985}). Even though our noise coefficient $\eps$ may be sufficiently small to justify to use a large deviations based approximation, we choose $\sqrt{2\eps}\nabla\phi$, with $\phi$ given by (\ref{mfet100}) as a control variate. The estimates of the MFET for standard MC and CoV with a badly chosen control variate  are shown in Figure \ref{fig:mfetOU} for $N=10$ and $M=100$, and they demonstrate that the subpotimal control variate still produces reasonable results at almost no additional computational cost.

\section{Discussion}

We have analysed the relationship between adaptive importance sampling (AIS) and control variate (CoV) on path space as two possible variance reduction methods for rare event simulation. Formulating AIS within a risk-sensitive control framework, we have shown that the latter is the risk-neutral limit of the former. Besides this conceptual link, it has been demonstrated that the risk-sensitive formulation of AIS allows for controlling several moments of the quantity of interest at once. This property is useful to avoid known pathologies of AIS that have been observed in connection with sampling problems that involve unbounded stopping times (e.g. zero-variance estimator with a.s.~infinite stopping time).  


Future work ought to address numerical approximations to the CoV scheme, especially non-asymptotic performance bounds when only a bad (i.e. suboptimal) control is available; cf.~\cite{roussel2019,vaes2023}. In practice, it is likely that a numerical approximation will be obtained from solving a simplified or lower-dimensional system, therefore it is crucial to analyse the robustness of the CoV estimator under suboptimal controls; the robustness issue is especially important if the length of the trajectories used for sampling is not bounded a priori. As a final remark, we emphasize that the (empirically) observed robustness of the CoV estimator under suboptimal controls is in stark contrast to the brittleness of high-dimensional AIS under suboptimal controls as has been pointed out by various authors, e.g. \cite{suboptimalIS,bengtsson2005curse,agapiou2015importance}. Related studies for CoV are still at its infancies, see e.g.~\cite{oates2023,belomestny2024} and references therein.

\begin{ack}
This work was partly supported by the DFG Collaborative Research Center 1114 “Scaling Cascades in Complex
Systems”, project No.235221301, Projects A05 ``Probing scales in equilibrated systems by optimal nonequilibrium forcing''.
One author (C.H.) acknowledges support by the German Federal Government, the Federal Ministry of Education and
Research and the State of Brandenburg within the framework of the joint project ``EIZ: Energy Innovation Center''
(project numbers 85056897 and 03SF0693A) with funds from the Structural Development Act (Strukturstärkungsge-
setz) for coal-mining regions
\end{ack}

%
%

\appendix

\section{Stochastic control problem}\label{sec:fbsde}

The sketch of proof of Theorem \ref{thm:is2cv} is based on a stochastic representation of the dynamic programming (or: HJB) equation associated with the optimal change of measure from $P$ to $Q^*$. For the sake of the argument, we consider only the exit problem:
\begin{ass}
        Let $B\subset\R^n$ be an open bounded set with smooth (at least $C^3$) boundary and $\tau=\tau_B$ its first exit time. We further assume: 
        \begin{itemize}
           \item $P(\tau<\infty)=1$, and 
           \[
           \bE\!\left[\int_0^\tau |u^*_t|^{2}dt\right]<\infty\ ,
           \]    
           \item $g=0$ and $f\in C^1(\overline{B},\R)$ is non-negative.
        \end{itemize}
\end{ass}

We let $V^\alpha(t,x)=\min_{u\in A}J(u;t,x)$ denote the value function associated with (\ref{cost}), where $A$ denotes the set of admissible Markovian controls, for which (\ref{controlledSDE}) has a unique strong solution strong; since the drift vector field $b$ in (\ref{controlledSDE}) is time-homogeneous, the strong Markov property of $X^u$ implies that $V^\alpha$ is independent of $t$, therefore we redefine 
\[V^\alpha(x):=V^\alpha(t,x)\,.\] 
By \cite[Thm. 5.1]{fleming2006}, $V^\alpha$ is the classical solution to the HJB equation 
\begin{equation}\label{hjb}
\begin{aligned}
        \cL V^\alpha - \frac{\alpha}{2}|\nabla V^\alpha|^2_{\sigma\sigma^T} + f & = 0 \,,\quad x\in B\\
        V^\alpha & = 0\,,\quad x\in\partial B\,,
\end{aligned}
\end{equation}
where $|h|_a=\sqrt{h\cdot ah}$ for some symmetric positive definite matrix $a\in\R^{n\times n}$ denotes the weighted Euclidean norm of a vector $h\in\R^n$. 

\subsection{Forward-backward SDE representation}

To study the convergence of the quantity of interest and the associated optimal control, it is convenient to consider the forward-backward SDE (in brief: FBSDE) representation of the dynamic programming equation. To this end, we define the two continuous, adapted processes
\begin{equation}
	Y_t^u = V^\alpha(X^u_t)\,,\quad Z^u_t = \sigma^T\nabla V^\alpha(X^u_t)\,,
\end{equation}
with $X^u$ solving the controlled SDE (\ref{controlledSDE}). We moreover recast (\ref{hjb}) as
\begin{equation*}
\begin{aligned}
        \cL^u V^\alpha - \frac{\alpha}{2}|\nabla V^\alpha|^2_{\sigma\sigma^T} - \sigma u \cdot \nabla V+ f & = 0 \,,\quad x\in B\\
        V^\alpha & = 0\,,\quad x\in\partial B\,.
\end{aligned}
\end{equation*}
Here $\cL^u$ denotes the infinitesimal generator associated with (\ref{controlledSDE}). Now, using It\^o's formula and the last equation,
\begin{equation*}
	\begin{aligned}
	Y^u_t  & = -\int_t^\tau (\cL^u V^\alpha)(X_s)ds - \int_t^\tau \sigma^T\nabla V(X^u_s) \cdot dW_s\\
	& = \int_t^\tau \left(f - \frac{\alpha}{2}|Z^u_s|^2 - u_s\cdot Z^u_s \right)ds - \int_t^\tau Z^u_s dW_s\,.
	\end{aligned}
\end{equation*}
In other words, the triple $(X^u,Y^u,Z^u)$ solves the FBSDE
\begin{equation}\label{fbsde}
	\begin{aligned}
			dX^u_t &= (b(X^u_t)+\sigma u_t)\,dt + \sigma dW_t\\
		dY^u_t & = \left(-f(X^u_t) + \frac{\alpha}{2}|Z^u_t|^2 + u_t Z^u_t\right)dt + Z^u_t dW_t\,,
			\end{aligned}
\end{equation}
that is equipped with the boundary conditions 
\begin{equation}\label{fbsde2}
    X^u_0=x\,,\quad  Y^u_\tau=0\,. 
\end{equation}

\subsection{Limit equation}

Equations (\ref{fbsde})--(\ref{fbsde2}) are the FBSDE representation of the CGF $\gamma$ under the controlled dynamics (i.e. under the measure $Q=Q(u)$.). Formally, substituting $u=u^*$ and setting $\alpha=0$, we obtain the FBSDE 
\begin{equation}\label{fbsdelim}
	\begin{aligned}
		dX_t &= b(X_t)\,dt + \sigma dW_t\,,\quad X_0=x\\
		dY_t & = -f(X_t)dt + Z_t dW_t\,, \quad Y_\tau=0\,,
	\end{aligned}
\end{equation}
where we have used that $u^*\to 0$ in $L^2$ as $\alpha\to 0$. By Dynkin's formula (\cite[Thm.~7.4.1]{oeksendal2003}), the FBSDE (\ref{fbsdelim}) represents the linear boundary value problem 
\begin{equation}\label{linBVP}
  \begin{aligned}
         -\cL\phi(x) & = f\,,\quad x\in B \\
         \phi(x) & = 0\,,\quad x\in \partial B\,,
    \end{aligned}
\end{equation}
the solution of which agrees with
\begin{equation}
    \phi(x) = \bE\!\left[\int_0^\tau f(X_t)\,dt\Big|X_0=x\right]
\end{equation}

\section{Control variate limit of AIS}

The FBSDE representation of the CGF is the basis for the proof of Theorem \ref{thm:is2cv}. We briefly sketch the main argument: Since the CGF agrees with the value function, the essential step is the uniform convergence of the (controlled) value function $V^\alpha$ to $\phi$ as $\alpha\to 0$ and the convergence $\nabla V^\alpha\to \nabla\phi$ of its gradients (which implies convergence of the AIS estimator to the control variate estimator).    

\subsection{Uniform convergence of the estimator}

By  \cite[Lemma A.2]{KNH18} the uniform convergence of the value function and its derivatives on any compact subset of $B$ is implied by the convergence $(Y^u,Z^u)\to(Y,Z)$ in $L^2$.  Since the FBSDE (\ref{fbsde}) associated with (\ref{hjb}) has a quadratic nonlinearity (i.e.~a non-Lipschitz right hand side), the trick is to first transform it to a linear FBSDE and then pass to the limit. We briefly sketch the key steps:  

\begin{enumerate}
    \item Setting $H_t=e^{-\alpha Y^u_t}$, It\^o's formula implies that 
    \begin{align*}
        dH^u_t & = -\alpha H^u_t dY^u_t + \frac{\alpha^2}{2}H^u_t |Z_t^u|^2\,dt\\
        & = a_t\,H^u_t\,dt + u_t K^u_t + K^u_t\,dW_t\,, 
    \end{align*}
    where we have used (\ref{fbsde}) and introduced the shorthands $a_t=\alpha f(X_t^u)$ and $K^u_t=-\alpha H^u_t Z_t^u$.
    \item The BSDE for the pair $(H,K)$ is linear, and it has the explicit solution \cite[Sec.~2.1.1]{bsde}  
    \[
    H^u_t=\bE\!\left[e^{-\alpha\int_t^\tau f(X_t^u)\,dt} L^{-1}_{t,\tau}\Big| X_t\right],
    \]
    where 
    \[
    L_{t,\tau} = \exp\left(\int_t^\tau u_t \cdot dW_t + \frac{1}{2}\int_t^\tau \left|u_t\right|^2 dt\right)\,.
    \]
    \item Here $\tau=\tau^u$ is the first exit time (FET) associated with the controlled process. Switching back to the reference measure $P$, generated by the uncontrolled dynamics (\ref{sde}), where $\tau=\tau^0$ now becomes again the uncontrolled FET, we have 
    \[
    H_t=\bE\!\left[e^{-\alpha\int_t^\tau f(X_t)\,dt} \Big| X_t\right],
    \]
    Note that we can drop the superscript $u$ because $H$ is independent of $u$. 
    \item Replacing now $\tau$ by the bounded stopping time $\tau_n=\min\{\tau,n\}$, $n\in\N$, and realising that $Y^u_t=\alpha^{-1}\log H^u_t$, with $H^u_t=H_t$ we can pass to the limit:
    \begin{align*}
        -\lim_{\alpha\searrow 0}\frac{1}{\alpha}\log H^u_t & = -\lim_{\alpha\searrow 0}\frac{1}{\alpha}\log \bE\!\left[e^{-\alpha\int_t^{\tau_n} f(X_t)\,dt}\Big| X_t\right]\\
        & = \bE\!\left[\int_t^{\tau_n} f(X_t)\,dt\Big| X_t\right]\,,
    \end{align*}
    where the last expression converges to $\phi(X_t)$ as $n\to\infty$ since $f$ is uniformly bounded on $\overline{B}$ which implies that the limit commutes with the expectation. 

    \item This implies that $Y^u_t\to Y_t$ uniformly in $[0,T]$ for arbitrary $T>0$ as $\alpha\to 0$. Then, by applying stability estimates \cite[Prop.~2.4]{Kobylanski2000} for BSDEs with quadratic nonlinearity, it can be shown that $Z^u\to Z$ in $L^2$, which yields convergence of AIS to the zero-variance control variate estimator.  
\end{enumerate}

\subsection{Zero-variance property}

The zero-variance property follows from the fact that $\phi$ solves the linear boundary value problem (\ref{linBVP}), since
\[
\phi(X_\tau) - \phi(x) = \int_0^\tau \cL\phi(X_t)\,dt + \int_0^\tau \sigma^T\nabla\phi(X_t)\cdot dW_t\,,
\]
which, by (\ref{linBVP}), implies almost surely that 
\[
\phi(x) = \int_0^\tau f(X_t)\,dt - \int_0^\tau \sigma^T\nabla\phi(X_t)\cdot dW_t\,.
\]
This concludes the sketch of proof of Theorem \ref{thm:is2cv}.

\end{document}